\documentclass [a4 paper,12pt]{article}
\usepackage{mathrsfs}

\usepackage[usenames]{color}
\usepackage{bbm,amsmath,amssymb,amsfonts,amsthm,graphics,graphicx}

\newtheorem{lem}{Lemma}
\newtheorem{thm}[lem]{Theorem}
\newtheorem{cor}[lem]{Corollary}

\DeclareMathOperator{\En}{\mathscr{E}}

\DeclareMathOperator{\si}{\sigma}
\DeclareMathOperator{\la}{\lambda}

\title{The asymptotic value of graph energy for random graphs with degree-based weights\footnote{Supported by NSFC No.11871034, 11531011 and NSFQH No.2017-ZJ-790.}}
\author{
\small  Xueliang Li$^{1,4}$, Yiyang Li$^2$, Jiarong Song$^3$\\
\small $^1$Center for Combinatorics and LPMC \\
\small Nankai University, Tianjin 300071, China \\
\small Email: lxl@nankai.edu.cn\\
\small $^2$International Institute, China Construction Bank\\
\small Xicheng District, Beijing 100033, China \\
\small Email: liyiyangnk@163.com\\
\small $^3$School of International Trade and Economics\\
\small University of International Business and Economics\\
\small Chaoyang District, Beijing 100029, China\\
\small Email: songjiarong9@ubie.edu.cn\\
\small $^4$School of Mathematics and Statistics\\
\small Qinghai Normal University\\
\small Xining, Qinghai 810008, China}
\date{}

\begin{document}
\maketitle

\begin{abstract}
In this paper, we investigate the energy of a weighted random graph $G_p(f)$ in $G_{n,p}(f)$, in which each edge $ij$ takes the weight $f(d_i,d_j)$, where $d_v$ is a random variable, the degree of vertex $v$ in the random graph $G_p$ of the Erd\"{o}s--R\'{e}nyi random graph model $G_{n,p}$, and $f$ is a symmetric real function on two variables. Suppose $|f(d_i,d_j)|\leq C n^m$ for some constants $C, m>0$, and $f((1+o(1))np,(1+o(1))np)=(1+o(1))f(np,np)$. Then, for almost all graphs $G_p(f)$ in $G_{n,p}(f)$, the energy of $G_p(f)$ is $(1+o(1))f(np,np)\frac{8}{3\pi}\sqrt{p(1-p)}\cdot n^{3/2},$ where $p\in(0,1)$ is any fixed and independent of $n$. Consequently, with this one basket we can get the asymptotic values of various kinds of graph energies of chemical use, such as Randi\'c energy, ABC energy, and energies of random matrices obtained from various kinds of degree-based chemical indices. \\[3mm]

\end{abstract}

\section{Introduction}

Throughout this paper, $G=(V,E)$ denotes a simple graph with vertex set $V$ and edge set $E$. Let $|V(G)|=n$. For a vertex $i\in V$, $d_i$ denotes the degree of $i$ in $G$, and if an edge $e\in E$ of $G$ has two end-vertices $i$ and $j$, we say that $e = ij \in E$. Let $A=(a_{ij})$ be the adjacency matrix of $G$. Assume that $f(x,y)$ be a symmetric real function in two variables $x$ and $y$. Consider the weighted graph $G(f)$ with edge weight $f(d_i,d_j)> 0$ for each edge $ij$ of $G$. The adjacency matrix of the weighted graph $G(f)$ is denoted by $A_G(f)$, or simply $A(f)$ if no confusion occurs. Since $f(x,y)$ is symmetric, $A(f)$ is a symmetric matrix. Suppose $\lambda_1(f),\ldots,\lambda_n(f)$ are the eigenvalues of $A(f)$. Then the energy of the weighted graph $G(f)$ is, as usual \cite{G}, defined as
$$
\En(G(f))=\En(A(f))=\sum_{i=1}^n|\lambda_i(f)|.
$$
More results on the energy of graphs can be found in \cite{glz, LSG}.

We will consider the energy of a random graph $G_p$ in the Erd\"{o}s--R\'{e}nyi random graph model $G_{n,p}$ \cite{BB}, in which the edges are taken independently with probability $p\in(0,1)$. Denote the corresponding adjacency matrix by $A(G_p)$. From a random graph $G_p$ in $G_{n,p}$, we can get a weighted random graph $G_p(f)$ in $G_{n,p}(f)$ with adjacency matrix $A(G_p(f))$, in which each random edge $ij$ has a weight $f(d_i,d_j)$. Since the degree $d_i$ of a vertex $i$ is now a random variable, the weight $f(d_i,d_j)$ of each edge is also a random variable. Because the weighted random graph $G_p(f)$ comes from the random graph $G_p$, $A(G_p(f))$ is a random matrix in which each $(i,j)$-entry $A(G_p(f))_{ij}$ takes weight $f(d_i,d_j)$ with probability $p\in(0,1)$.

Throughout this paper we always assume that the probability value $p\in (0,1)$ is any fixed and independent of $n$. In the following, we will focus on the property of the eigenvalues of the random matrix $A(G_p(f))$.

Relying on the moments method, Wigner \cite{W55,W58} considered the limiting spectral distribution (LSD for short) of a kind of random matrix and $\bar{A}(G_p)=A(G_p)-p(J-I)$ is a special case, where $I$ is the
unit matrix  and $J$ is the matrix in which all entries equal 1.The LSD is called semi-circle law.

\begin{thm}\label{Thm-1} \cite{W58}
Let $\Phi_{\mathbf{\bar A(G_p)}}(x)=
\frac{1}{n}\cdot\#\{\la_{i}\mid\la_{i}\le x, ~i=1,2,\ldots, n\}$(ESD).
Then
$$
\lim_{n\rightarrow\infty}\Phi_{n^{-1/2}\mathbf{\bar A(G_p)}}(x)=\Phi(x)
\mbox{ a.s. }
$$
{\it i.e.,} with probability 1, $\Phi_{n^{-1/2}\mathbf{\bar A(G_p)}}(x)$ converges weakly to a distribution $\Phi(x)$ as $n$ tends to infinity.
$\Phi(x)$ has the density
$$
\phi(x)=
\frac{1}{2\pi\si_2^2}\sqrt{4\si_2^2- x^2 }~\mathbf{1}_{|x|\le2\si_2},
$$
where $\si_2=\sqrt{p(1-p)}$.
\end{thm}

Based on this theorem, the following result was obtained in \cite{dll, N}

\begin{thm}\label{Thm-D} \cite{dll, N}
For almost all graphs $G_p$ in $G_{n,p}$, the energy of $G_p$ is
$$(1+o(1))\frac{8}{3\pi}\sqrt{p(1-p)}n^{3/2},$$
where $p\in (0,1)$ is any fixed and independent of $n$.
\end{thm}

Of course, one can calculate the exact value of energy for each graph, and different graphs may have different values of energy. But from the probability point of view, almost all graphs have the same value of energy. Similarly, in this paper we shall get the following result relying on Winger's methods.

\begin{thm}\label{Thm-F}
Let $f(x,y)$ be a symmetric real function such that
$|f(d_i, d_j)|, 1/|f(d_i, d_j)|\leq C n^m$ for some constants $C, m>0$, and $f((1+o(1))np,(1+o(1))np)=(1+o(1))f(np,np)$
where $p\in (0,1)$ is any fixed and independent of $n$ and $f(np,np)\neq 0$. Then
for almost all weighted random graphs $G$ in $G_{n,p}(f)$, the energy of
$G$ is
$$
\En(A(G_p(f)))=(1+o(1))f(np,np)\frac{8}{3\pi}\sqrt{p(1-p)}\cdot n^{3/2}.
$$
\end{thm}

We can see that even if each edge is given a weight by a function of
the degrees of the two end-vertices, the values of the energy for most weighted random graphs are intensely concentrated in a small interval, and the values are mainly determined by $f(np,np)$. Applying this result, with one basket we can get the energies of almost all random matrices defined by various kinds of degree-based chemical indices.

\section{The asymptotic value of energy for random graphs $G_p(f)$ with degree-based weights}

In this section, we will give a proof of Theorem \ref{Thm-F}, by estimating the asymptotic value of energy for the random adjacency matrices $A(G_p(f))$ of weighted random graphs $G_p(f)$ in $G_{n,p}(f)$ with a weight function $f$ on degree-based variables.

In fact, the research on the spectral distributions of random
matrices is rather abundant and active, which can be dated back to Wishart \cite{Wi}. We refer the readers to \cite{B, De, fk, Me, V} for an overview and some spectacular progress in this field. One important
achievement in that field is the Wigner's semi-circle law \cite{W58}. Moreover, the spectral theory of random matrices has revealed deep links with many fields of mathematics, theoretical physics and chemistry.
In these researches, the key point is to computing the moments of a random matrix $A$.

Let
$$M_k=\int x^k d \Phi_{\mathbf{ A}}(x).$$
Then
$$\textbf{E}M_k=\textbf{E}\frac{1}{n}\displaystyle \sum_{i=1}^n\lambda_i^k=\textbf{E}\frac{1}{n} Tr(A^k).$$

Our proof of Theorem \ref{Thm-F} goes as follows.

{\bf Proof of Theorem \ref{Thm-F}:}
To get the energy of $A(G_p(f))$, we also use the method of computing moments, and study the mathematical expectations.
In the following, we shall try to establish the relation between $M_k(\widetilde{\bar A(G_p(f)))}$ and $M_k(\bar A(G_p))$. Then, by the moment method, we can get the eigenvalue distribution of $\widetilde{\bar A(G_p(f))}$.

Consider a new matrix which is defined as
$$\widetilde{\bar A(G_p(f))}=A(G_p(f))/f(np,np)-p(J-I).$$
Denote the element in $\widetilde{\bar A(G_p(f))}$ by $\tilde f(d_i,d_j)$ and it equals to $f(d_i,d_j)/f(np,np)-p$ with probability $p$ or $-p$ with probability $1-p$. We show that for any positive integer number $k$
\begin{equation}\label{eq0}
\textbf{E}M_k(\widetilde{\bar A(G_p(f))})\sim \textbf{E}M_k\bar A(G_p))  \mbox { a.s..}
\end{equation}

At first, we get the estimation of $\textbf{E}M_k(A(G_p(f)))$ as follows:
$$\begin{array}{lll}
  \textbf{E}M_k(\widetilde{\bar A(G_p(f))})&=&\textbf{E} \frac{1}{n} Tr(\widetilde{\bar A(G_p(f))})^k \\
   &=& \frac{1}{n} \displaystyle\sum_{i_1=1}^n\displaystyle\sum_{i_2=1}^n\cdots\displaystyle\sum_{i_k=1}^n \textbf{E}(\tilde f(d_{i_1},d_{i_2})\cdot \tilde f(d_{i_2},d_{i_3})\cdots \tilde f(d_{i_k},d_{i_1})),
 \end{array}$$
where $\{i_1i_2\ldots i_ki_1\}$ is a  $k$-closed walk. All these $k$-closed walks form a set, denoted by $W_k$. Each closed walk provides
a contribution
$$\textbf{E}(\tilde f(d_{i_1},d_{i_2})\cdot \tilde f(d_{i_2},d_{i_3})\cdots \tilde f(d_{i_k},d_{i_1}))$$
to the summation.

If the weights of the edges are independent, it was shown in \cite{ksm} that the eigenvalue distribution converges to a limiting measure as $n$ tends to infinity by the moments and the resolvent techniques. And the limiting measure is much more complicated than semi-circle law.
In our case, the random variable $\tilde f(d_i,d_j)$ is related to the adjacency of vertices $i$ and $j$. We could not use the conclusions in \cite{ksm}.

Now, we focus on the estimation of $\textbf{E}M_k$. Recall that $W_k$ is the set of closed walks of $k$ steps over the set $\{1,\cdots, n\}$. Hence,
\begin{equation}\label{eq1}
\textbf{E}M_k(\widetilde{\bar A(G_p(f))})=\frac{1}{n}\sum_{w\in W_k} \textbf{E}(\tilde f(d_{i_1},d_{i_2})\cdot \tilde f(d_{i_2},d_{i_3})\cdots \tilde f(d_{i_k},d_{i_1})).
\end{equation}
Set $a_{i_li_{l+1}}=1-p$ with probability $p$ and $a_{i_li_{l+1}}=-p$ with probability $1-p$.
Then,
$$\begin{array}{lll}
  \textbf{E}M_k(\widetilde{\bar A(G_p(f))})&=&\frac{1}{n}\displaystyle\sum_{w\in W_k}\displaystyle\sum (\tilde f(d_{i_1},d_{i_2})\cdot \tilde f(d_{i_2},d_{i_3})\cdots \tilde f(d_{i_k},d_{i_1}))\\
  & &\cdot \textbf{p}(d_{i_1},d_{i_2},\cdots,d_{i_k},\{i_1i_2\cdots i_ki_1\}\subseteq s(w))\\
    &=&\frac{1}{n}\displaystyle\sum_{w\in W_k}\displaystyle\sum (\tilde f(d_{i_1},d_{i_2})\cdot \tilde f(d_{i_2},d_{i_3})\cdots \tilde f(d_{i_k},d_{i_1}))\\
  & &\cdot \textbf{p}(d_{i_1},d_{i_2},\cdots,d_{i_k}|\{i_1i_2\cdots i_ki_1\}\subseteq s(w))\\
  & &\cdot\textbf{p}(\{i_1i_2\cdots i_ki_1\}\subseteq s(w))\\
  &=&\frac{1}{n}\displaystyle\sum_{w\in W_k} \textbf{p}(\{i_1i_2\cdots i_ki_1\}\subseteq s(w))\cdot a_{i_1i_2}\cdot a_{i_2i_3}\cdots a_{i_ki_1}\\
  & &\displaystyle\sum \{\frac{\tilde f(d_{i_1},d_{i_2})}{a_{i_1i_2}}\cdot \frac{\tilde f(d_{i_2},d_{i_3})}{a_{i_2i_3}}\cdots \frac{\tilde f(d_{i_k},d_{i_1}))}{a_{i_ki_1}}\cdot \\
  & &\textbf{p}(d_{i_1},d_{i_2},\cdots,d_{i_k}|\{i_1i_2\cdots i_ki_1\}\subseteq s(w)\},
 \end{array}$$
where $s(w)$ denotes the edge set of a closed walk $w$ taken from $G_p$ in the random graph model $G_{n,p}$. And we assume that the edges of $w$ form a simple subgraph $T_w=(v(w),s(w))$.
We should notice that $|v(w)|,|s(w)|\leq k$.

Recall that in the Erd\"{o}s--R\'{e}nyi model, the degree of any vertex $d_i$ is subject to binomial distribution $B(n,p)$. For large $n$,
$$
d_i\sim N(np,np(1-p)), {\text{ and } }  \textbf{p}(d_i)=\frac{1}{\sqrt{2\pi np(1-p)}}\cdot\exp^{-\frac{(d_i-np)^2}{2np(1-p)}},
$$
and $d_i=(1+o(1))np$ with probability 1.

Remember that we have assumed $|f(d_i,d_j)|,1/|f(d_i, d_j)|\leq Cn^m$, thus $|\tilde f(d_i,d_j)|<C^2n^{2m}$. For any vertex $i_l\in v(w)$ in each closed walk $w$,
$$
\begin{array}{lll}
& &\textbf{p}(d_{i_1},d_{i_2},\cdots,d_{i_k}|\{i_1i_2\cdots i_ki_1\}\subseteq s(w))\\
& &= \frac{\textbf{p}(d_{i_1},d_{i_2},\cdots,d_{i_k},\{i_1i_2\cdots i_ki_1\}\subseteq s(w))}{\textbf{p}(a_{i_1i_2},a_{i_2i_3},\cdots,a_{i_ki_1})}\\
& &\leq \frac {\textbf{p}(d_{i_l})}{p^{|s(w)|}(1-p)^{|s(w)|}},
\end{array}
$$

$$
\begin{array}{lll}
& &\big |\sum\cdots\sum_{d_{i_1},d_{i_2},\cdots,d_{i_{|v(w)|}},d_{i_l}<{np-n^{\frac{3}{4}}}} \frac{\tilde f(d_{i_1},d_{i_2})}{a_{i_1i_2}}\cdots \frac{\tilde f(d_{i_k},d_{i_1})}{a_{i_ki_1}}\\
& &\cdot \textbf{p}(d_{i_1},d_{i_2},\cdots,d_{i_k}|\{i_1i_2\cdots i_ki_1\}\subseteq s(w))\big |\\
&<&\sum\cdots\sum_{d_{i_1},d_{i_2},\cdots,d_{i_{|v_w|}},d_{i_l}<{np-n^{\frac{3}{4}}}} \big |\frac{\tilde f(d_{i_1},d_{i_2})}{a_{i_1i_2}}\cdots \frac{\tilde f(d_{i_k},d_{i_1})}{a_{i_ki_1}}\big |\cdot \textbf{p}(d_{i_l})/({p^{|s(w)|}}(1-p)^{|s(w)|})\\
&<& C^{2k}n^{2mk}\cdot \sum\cdots\sum_{d_{i_1},d_{i_2},\cdots,d_{i_{|v(w)|}},d_{i_l}<{np-n^{\frac{3}{4}}}} \frac{1}{\sqrt{2\pi np(1-p)}}\cdot\exp^{-\frac{n^{\frac{3}{2}}}{2np(1-p)}}/({p^{2|s(w)|}(1-p)^{2|s(w)|}})\\
&<& C^{2k}n^{2mk}\cdot n^k \cdot \exp^{-2{\sqrt{n}}}/{(p^{2|s(w)|}(1-p)^{2|s(w)|})}\\
&\rightarrow& 0 \mbox{ as } n\rightarrow\infty.
\end{array}
$$
Since $p(1-p)\leq 1/4$, the last inequality holds if $n$ is large enough.

Similarly,
$$
\begin{array}{lll}
\sum\cdots\sum_{d_{i_1},d_{i_2},\cdots,d_{i_{|v(w)|}},d_{i_l}>{np+n^{\frac{3}{4}}}} \frac{\tilde f(d_{i_1},d_{i_2})}{a_{i_1i_2}}\cdots \frac{\tilde f(d_{i_k},d_{i_1})}{a_{i_ki_1}}\\
\cdot \textbf{p}(d_{i_1},d_{i_2},\cdots,d_{i_k}|\{i_1i_2\cdots i_ki_1\}\subseteq s(w)))
\rightarrow 0 \mbox{ as } n\rightarrow\infty.
\end{array}
$$

Recall that $$\mathbf{P}(S_1\cap S_2\cap\cdots \cap S_t)\geq 1-(\mathbf{P}(\bar{S_1})+\mathbf{P}(\bar{S_2})+\cdots +\mathbf{P}(\bar{S_t}))$$
for any probability measure $\mathbf{P}$ and measurable sets $S_1,\cdots, S_t$.

Then, relying on the distribution of the vertex degree $d_i$, one can see that
$$
\begin{array}{lll}
& &\displaystyle\sum\cdots\sum_{np-n^{\frac{3}{4}}\leq d_{i_1},d_{i_2},\cdots,d_{i_{|v(w)|}}\leq{np+n^{\frac{3}{4}}}}\textbf{p}(d_{i_1},d_{i_2},\cdots, d_{i_k}|\{i_1i_2\cdots i_ki_1\}\subseteq s(w))\\
& &\geq 1-\displaystyle\sum_{l=1}^k\displaystyle\sum_{|d_{i_l}-np|\geq n^{\frac{3}{4}}}\textbf{p}(d_{i_1}|\{i_1i_2\cdots i_ki_1\}\subseteq s(w))\\
& &\geq 1-k\cdot\displaystyle\sum_{|d_{i_l}-np|\geq n^{\frac{3}{4}}} \frac{1}{\sqrt{2\pi np(1-p)}}\cdot\exp^{-\frac{n^{\frac{3}{2}}}{2np(1-p)}}/{p^{|s(w)|}}\\
& &=1 \mbox{,   as } n\rightarrow\infty.
\end{array}
$$
Hence, for each closed walk $w\in W_k$,
$$
\begin{array}{lll}
& &\displaystyle\sum\cdots\sum_{np-n^{\frac{3}{4}}\leq d_{i_1},d_{i_2},\cdots,d_{i_{|v(w)|}}\leq{np+n^{\frac{3}{4}}}} \frac{\tilde f(d_{i_1},d_{i_2})}{a_{i_1i_2}}\cdots \frac{\tilde f(d_{i_k},d_{i_1})}{a_{i_ki_1}}\\
& &\cdot \textbf{p}(d_{i_1},d_{i_2},\cdots,d_{i_k}|\{i_1i_2\cdots i_ki_1\}\subseteq s(w))\\
&=&\sum\cdots\sum_{np-n^{\frac{3}{4}}\leq d_{i_1},d_{i_2},\cdots,d_{i_{|v(w)|}}\leq{np+n^{\frac{3}{4}}}}(1+o(1))\\
& &\textbf{p}(d_{i_1},d_{i_2},\cdots, d_{i_k}|\{i_1i_2\cdots i_ki_1\}\subseteq s(w))\\
&=&\displaystyle 1  \mbox{,   as } n\rightarrow\infty.
\end{array}
$$
Therefore,
$$
\begin{array}{lll}
  \textbf{E}M_k(\widetilde{\bar A(G_p(f))})&\sim &\frac{1}{n}\displaystyle\sum_{w\in W_k} \textbf{p}(\{i_1i_2\cdots i_ki_1\}\subseteq s(w))\cdot a_{i_1i_2}\cdot a_{i_2i_3}\cdots a_{i_ki_1}\\
  &=&\frac{1}{n}\sum_{w\in W_k} \textbf{E}(a_{i_1i_2}\cdot a_{i_2i_3}\cdots a_{i_ki_1})\\
  &=&\textbf{E}M_k(\bar A(G_p))\mbox { a.s.}
\end{array}
$$
That means that $\lim_{n\rightarrow\infty} \textbf{E}M_k(n^{-1/2}\widetilde{\bar A(G_p(f))})=\lim_{n\rightarrow\infty} \textbf{E} M_k(n^{-1/2}\bar A(G_p)) \mbox { a.s}$.

To meet the requirements of moment approach referred in \cite{B}, we calculate the estimation of $\textbf{Var}(M_k(\widetilde{\bar A(G_p(f))})$ in the following.
It is well known that
$$
\begin{array}{lll}
  \textbf{Var}M_k(\widetilde{\bar A(G_p(f))})&=&\textbf{E}M_k^2(\widetilde{\bar A(G_p(f))})-\textbf{E}{M_k(\widetilde{\bar A(G_p(f))})}^2.
\end{array}
$$
Let us focus on $\textbf{E}M_k^2(\widetilde{\bar A(G_p(f))})$:
$$
\begin{array}{lll}
  \textbf{E}M_k^2(\widetilde{\bar A(G_p(f))})&=& \textbf{E}\frac{1}{n^2}\displaystyle\lbrack\sum_{w\in W_k}(\tilde f(d_{i_1},d_{i_2})\cdot \tilde f(d_{i_2},d_{i_3})\cdots \tilde f(d_{i_k},d_{i_1}))\rbrack ^2 \\
  &=&\frac{1}{n^2} \displaystyle\sum\sum_{w_i,w_j\in W_k}\textbf{E}(\tilde f(d_{i_1},d_{i_2})\cdot \tilde f(d_{i_2},d_{i_3})\cdots \tilde f(d_{i_k},d_{i_1})\\
  & &\cdot \tilde f(d_{j_1},d_{j_2})\cdot \tilde f(d_{j_2},d_{j_3})\cdots \tilde f(d_{j_k},d_{j_1})),
\end{array}
$$
Repeating the process of estimation of $\textbf{E}M_k(\widetilde{\bar A(G_p(f))})$, we can get
$$
\textbf{E}M_k^2(\widetilde{\bar A(G_p(f))})=(1+o(1))\frac{1}{n^2}\displaystyle\sum\sum_{w_i,w_j\in W_k} \textbf{E}(a_{i_1i_2}\cdots a_{i_ki_1}\cdot a_{j_1j_2}\cdots a_{j_kj_1}).
$$
Thus,
\begin{equation}\label{eq2}
\textbf{Var}(M_k(\widetilde{\bar A(G_p(f))})=(1+o(1)\textbf{Var}(M_k(\bar A(G_p))) \mbox{ a.s.}
\end{equation}
Therefore, the limiting distribution of the eigenvalues of $n^{-1/2}\widetilde{\bar A(G_p(f))}$ is the same as $n^{-1/2}\bar A(G_p)$. We refer the readers to \cite{B, BDJ} for more information on the moment method. As in \cite{dll}, by Theorem \ref{Thm-1}, we can get that $\sum_{i=1}^n |\lambda_i(n^{-1/2}\widetilde{\bar A(G_p(f))})|\sim \int |x|d \Phi(x)\cdot n=\frac{8}{3\pi}\sqrt{p(1-p)}\cdot n$ and $$\En(\widetilde{\bar A(G_p(f))}=(1+o(1))\frac{8}{3\pi}\sqrt{p(1-p)}\cdot n^{3/2} \mbox{ a.s.}$$

Let $\mathbf{X},\mathbf Y$ be real symmetric matrices of order $n$. Ky Fan's inequality \cite{KF} shows that $\sum_{i=1}^n |\lambda_i(\mathbf X+\mathbf Y)|\leq \sum_{i=1}^n |\lambda_i(\mathbf X)|+\sum_{i=1}^n |\lambda_i(\mathbf Y)|$.
Recall that
$$\widetilde{\bar A(G_p(f))}=A(G_p(f))/f(np,np)-p(J-I).$$
Consequently,
$$\En(\widetilde{\bar A(G_p(f))})-\En(-p(J-I))\leq\En(A(G_p(f))/f(np,np))\leq \En(\widetilde{\bar A(G_p(f))})+\En(p(J-I)).$$
It is easy to see that $\En(p(J-I))=\En(-p(J-I))=2p(n-1)$. Then, it follows that
$$
\begin{array}{lll}
  \En(A(G_p(f)))=(1+o(1))f(np,np)\frac{8}{3\pi}\sqrt{p(1-p)}\cdot n^{3/2}, \mbox{ a.s.}
 \end{array}
$$
which completes the proof of Theorem \ref{Thm-F}. \ \ \ \ \ \ \ \ \ \ \ \ $\blacksquare$
\\

Let $\lVert f\rVert=sup_x|f(x)|$. From \cite{B}, the ESDs of two symmetric matrices $\mathbf{X},\mathbf{Y}$ satisfy
$$\lVert \Phi_\mathbf{X}(x)-\Phi_\mathbf{Y}(x)\rVert \leq \frac{rank(\mathbf{X}-\mathbf{Y})}{n}.$$
Thus, we have
$$\lVert \Phi_{\frac{1}{n^{1/2}f(np,np)}{A(G_p(f))}}(x)-\Phi_{\frac{1}{n^{1/2}f(np,np)}{A(G_p(f))}-\frac{1}{n^{1/2}}pJ}(x)\rVert \leq \frac{rank(\frac{1}{n^{1/2}}pJ)}{n}=\frac{1}{n}.$$
Since $\lambda_1(I)=\cdots =\lambda_n(I)=1$, one can see from the definition of ESD that
$$\Phi_{\frac{1}{n^{1/2}f(np,np)}{A(G_p(f))}-\frac{1}{n^{1/2}}pJ}(x)=\Phi_{\widetilde{\bar A(G_p(f))}}(x+p/\sqrt{n}).$$
As $n\rightarrow \infty$, $\lim_{n\rightarrow\infty}\Phi_{\frac{1}{n^{1/2}f(np,np)}{A(G_p(f))}}(x)=\lim_{n\rightarrow\infty}\Phi_{\widetilde{\bar A(G_p(f))}}(x).$
\begin{thm}\label{Thm-FF}
Let $A(G_p(f))$ be the adjacency matrix of a weighted random graph $G_p(f)$ in the weighted random graph model $G_{n,p}(f)$ with weight function $f(x,y)$ based on the vertex degrees of the random graph $G_p$ in the Erd\"{o}s--R\'{e}nyi random graph model $G_{n,p}$. Then the eigenvalue distribution of $A(G_p(f))$ enjoys
$$
\lim_{n\rightarrow\infty}\Phi_{\frac{1}{n^{1/2}f(np,np)}{A(G_p(f))}}(x)=\Phi(x)
\mbox{ a.s.,}
$$
where $|f(d_i,d_j)|,1/|f(d_i,d_j)|\leq C n^m$ for some constants $C, m>0$, and $f((1+o(1))np,(1+o(1))np)=(1+o(1))f(np,np)$ in which
$p\in (0,1)$ is any fixed number and independent of $n$.
\end{thm}

\noindent{\bf Remark 1:} Although we supposed $|f(d_i,d_j)|,1/|f(d_i,d_j)|\leq C n^m$   for some constants $C,m>0$, and $f((1+o(1))np,(1+o(1))np)=(1+o(1))f(np,np)$, the restrictions are not so strict. For example, these requirements hold for all polynomial functions. And in the next section, we can see that most degree-based matrices of chemical use meet the above restrictions. If $f$ is differentiable, then $\max_{-n^{3/4} \leq x\leq n^{3/4}}f'_x(np+x,np)
= o(n^{-3/4}f(np,np))$ is a necessary condition. Since $f$ is a symmetric function, $f'_x(np+x,np)=f'_x(np,np+x)$. By the mean value theorem, we can get that $f(np+x,np+y)=(1+o(1))f(np,np)$ for any $x,y\in [-n^{3/4}, n^{3/4}]$. We should notice that the differentiability of $f$ is not a
necessary condition for $f((1+o(1))np,(1+o(1))np)=(1+o(1))f(np,np)$.

\section{Application to the energy of $A(G_p(f))$ for degree-based weight functions $f$ of chemical use}

Topological indices have been proved to be very useful in various chemical disciplines, which are associated with chemical structures.
Degree-based indices are the kinds of often used indices. The general form is
$$
T(G)=\sum_{ij\in E(G)}f(d_i,d_j),
$$
where $T(G)$ denotes the topological index value of a graph $G$. In \cite{FGD}, the authors collected many kinds of topological indices, and each kind corresponds to a different weight function $f(d_i,d_j)$ on degree-based variables. Then, the degree-based weighted matrices can be obtained from various degree-based chemical indices of graphs. One can see that the function $f$ for each of these topological indices meets the requirements of Theorem \ref{Thm-F}. Hence, in this section, we get the asymptotic values of the energy for corresponding weighted random matrices defined by degree-based indices, directly from Theorem \ref{Thm-F}.

Here we would like to point out that topological indices in chemistry are used to represent structural properties of molecular graphs. Each index maps a molecular graph into a single number,
obtained by summing up all the edge-weights of a molecular graph. If we use a weighted adjacency matrix to represent the structure of a molecular graph with weights separately on edges, it will keep the
complete structural information of the graph. That is, a (2-dimensional) matrix is much more effective than a (0-dimensional) index for representing structural information. Even if we consider the
spectrum of a matrix, it gives us a (1-dimensional) information, better than a 0-dimensional information. So, to further study the algebraic properties of
these structural matrices should be a necessary job to do.

The first degree-based matrix perhaps could be the Zagreb matrix, coming from the so-called {\it Zagreb indices} $M_1$ and $M_2$. We denote them by $M_1(G)$ and $M_2(G)$, which are defined as follows.\\

(1) For $M_1(G)=(m_{ij})$, $m_{ij}=d_i+d_j$, if $ij\in E$; 0 otherwise.\\

(2) For $M_2(G)=(m_{ij})$, $m_{ij}=d_id_j$, if $ij\in E$; 0 otherwise.\\
\begin{cor}\label{cor1}
For almost all graphs $G$ in the Erd\"{o}s--R\'{e}nyi model $G_{n,p}$ with weights from the Zagreb indices,
$$\begin{array}{lll}
  \displaystyle\En(M_1(G))&=&(1+o(1))p\cdot \frac{16}{3\pi}\sqrt{p(1-p)}\cdot n^{5/2},\\

  \displaystyle\En(M_2(G))&=&(1+o(1))p^2\cdot \frac{8}{3\pi}\sqrt{p(1-p)}\cdot n^{7/2}.
 \end{array}$$

\end{cor}
It is easily seen that the maximum is achieved when $p=3/4$ and  $p=5/6$, respectively.\\

Another earlier index is the so-called {\it Randi\'c index} $R$. From it we can get the Randi\'c matrix $R(G)$ as follows. \\

(3) For $R(G)=(r_{ij})$, $r_{ij}= \frac {1} {\sqrt{d_id_j}}$, if $ij\in E$; 0 otherwise.\\
\begin{cor}\label{cor2}
For almost all graphs $G$ in the Erd\"{o}s--R\'{e}nyi model $G_{n,p}$ with weight from the Randi\'c index,
$$\begin{array}{lll}
  \En(R(G))&=&(1+o(1))\cdot \frac{8}{3\pi}\displaystyle\sqrt{\frac{1-p}{p}}\cdot n^{1/2}.\\
  \end{array}$$
\end{cor}
It is easily seen that when $p$ is small, the energy is large, but when $p$ is large the energy is small. This coincides with
the fact that complete graph has the minimum Randi\'c energy, and path has the maximum Randi\'c energy, conjectured in \cite{gutman}.\\

Later, Bollobas and Erd\'os introduced the concept of {\it general Randi\'c index} $R_\alpha$, and from which we can get a
general Randi\'c matrix as follows.\\

(4) For $R_\alpha(G)=(r_{ij})$, $r_{ij}=(d_id_j)^\alpha$, if $ij\in E$; 0 otherwise.\\
\begin{cor}\label{cor3}
For almost all graphs $G$ in the Erd\"{o}s--R\'{e}nyi model $G_{n,p}$ with weight from the general Randi\'c index $R_\alpha(G)$,
$$\begin{array}{lll}
  \En(R_\alpha(G))&=&(1+o(1))p^{2\alpha}\cdot \frac{8}{3\pi}\displaystyle\sqrt{p(1-p)}\cdot n^{3/2+2\alpha}.\\
  \end{array}$$
\end{cor}

It is easy to see that from (4) we can get (2) and (3) by setting $\alpha=1$ and $-\frac 1 2$, respectively.
It is easily seen that The maximum is achieved when $p=\frac {4\alpha+1} {4\alpha +2}$, where $\alpha\geq -1/4$ or $<-1/2$.
What about when $-1/2<\alpha < -1/4$ ? Actually, it is already not obvious when $\alpha <0$.\\

The {\it ABC-index} is another kind, from which we can get the ABC-matrix of a graph $G$ as follows.\\

(5) For $ABC(G)=(ABC_{ij})$, $ABC_{ij}=\frac {\sqrt{d_i+d_j-2}} {\sqrt{d_id_j}}$, if $ij\in E$; 0 otherwise.\\
\begin{cor}\label{cor4}
For almost all graphs $G$ in the Erd\"{o}s--R\'{e}nyi model $G_{n,p}$  with weight from the ABC-index $ABC(G)$,
$$\begin{array}{lll}
  \En(ABC(G))&=&(1+o(1))\cdot \frac{8\sqrt{2}}{3\pi}\sqrt{1-p}\cdot n.\\
  \end{array}$$
\end{cor}
It is easily seen that when $p$ is small, the energy is large, but when $p$ is large the energy is small. \\

Recently, Furtula introduced the {\it augmented Zagreb index (AZI)}, from which we can get the AZI-matrix as follows.\\

(6) For $AZI(G)=(AZI_{ij})$, $AZI_{ij}= (\frac {d_id_j} {d_i+d_j-2})^3$, if $ij\in E$; 0 otherwise.\\
\begin{cor}\label{cor5}
For almost all graphs $G$ in the Erd\"{o}s--R\'{e}nyi model $G_{n,p}$ with weight from the Augmented Zagreb index $AZI(G)$,
$$\begin{array}{lll}
  \En(AZI(G))&=&(1+o(1))\cdot \frac{p^3\sqrt{p(1-p)}}{3\pi}\cdot n^{9/2}.\\
  \end{array}$$
\end{cor}
It is easily seen that the maximum is achieved when $p=7/8$. \\

From the {\it Arithmetic-Geometric ($AG_1$) index}, we get the  arithmetic-geometric matrix $AG_1$ as follows.\\

(7) For $AG_1(G)=(AG_{ij})$, $AG_{ij}=\frac {2\sqrt{d_id_j}}{d_i+d_j}$, if $ij\in E$; 0 otherwise.\\
\begin{cor}\label{cor6}
For almost all graphs $G$ in the Erd\"{o}s--R\'{e}nyi model $G_{n,p}$ with weight from the arithmetic-geometric index $AG_1(G)$,
$$\begin{array}{lll}
  \En(AG_1(G))&=&(1+o(1))\frac{8}{3\pi}\sqrt{p(1-p)}\cdot n^{3/2}.\\
  \end{array}$$
\end{cor}
It is easily seen that the maximum is achieved when $p=1/2$.\\

From the {\it Harmonic index} $HI$, the {\it Harmonic} matrix is in the form:\\

(8) For $HI(G)=(HI_{ij})$, $HI_{ij}=\frac {2}{d_i+d_j}$, if $ij\in E$; 0 otherwise.\\
\begin{cor}\label{cor7}
For almost all graphs $G$ in the Erd\"{o}s--R\'{e}nyi model $G_{n,p}$ with weight from the {\it Harmonic} index $HI(G)$,
$$\begin{array}{lll}
  \En(HI(G))&=&(1+o(1))\frac{8}{3\pi}\sqrt{\frac{1-p}{p}}\cdot n^{1/2}.\\
  \end{array}$$
\end{cor}
It is easily seen that the behavior likes the Randi\'c energy.\\

From the {\it Sum-connectivity index} $SCI$, the {\it Sum-connectivity} matrix $SCI(G)$ is defined as\\

(9) For $SCI(G)=(SCI_{ij})$, $SCI_{ij}=\frac {1}{\sqrt{d_i+d_j}}$, if $ij\in E$; 0 otherwise.\\
\begin{cor}\label{cor8}
For almost all graphs $G$ in the Erd\"{o}s--R\'{e}nyi model $G_{n,p}$ with weight from the {\it Sum-connectivity} index $SCI(G)$,
$$\begin{array}{lll}
  \En(SCI(G))&=&(1+o(1))\frac{2\sqrt{2}}{3\pi}\sqrt{1-p}\cdot n.\\
  \end{array}$$
\end{cor}
It is easily seen that the behavior likes the ABC energy.\\

From the {\it First (modified) multiple Zagreb index} $\ln \prod_1$($\ln \prod_1^*$), we get the {\it First (modified) multiple Zagreb} matrix as follows.\\

(10a) For $\ln \prod_1(G)=(\ln {\prod_1}_{ij})$, $\ln {\prod_1}_{ij}=\frac{\ln d_i}{d_i}+\frac{\ln d_j}{d_j}$, if $ij\in E$; 0 otherwise.\\

(10b) For $\ln \prod_1^*(G)=(\ln {\prod_1^*}_{ij})$, $\ln {\prod_1^*}_{ij}=\ln (d_i+d_j)$, if $ij\in E$; 0 otherwise.
\begin{cor}\label{cor9}
For almost all graphs $G$ in the Erd\"{o}s--R\'{e}nyi model $G_{n,p}$ with weight from the {\it First (modified) multiple Zagreb} indices $\ln \prod_1(G)$ and $ \ln \prod_1^*(G)$,
$$\begin{array}{lll}
  \En(\ln \prod_1(G))&=&(1+o(1))\frac{16}{3\pi}\sqrt{\frac{1-p}{p}}\cdot n^{1/2}\ln n,\\
  \En(\ln \prod_1^*(G))&=&(1+o(1))\frac{8}{3\pi}\sqrt{p(1-p)}\cdot n^{3/2}\ln n.
  \end{array}$$
\end{cor}
It is easily seen that the behavior likes the Randi\'c energy and the adjacency energy ($p=1/2$), respectively.\\

From the {\it Second multiple Zagreb index} $\ln \prod_1$, we get the {\it Second multiple Zagreb} matrix as follows.\\

(11) For $\ln \prod_2(G)=(\ln {\prod_1}{ij})$, $\ln {\prod_2}{ij}=\ln d_i+\ln d_j$, if $ij\in E$; 0 otherwise.\\
\begin{cor}\label{cor10}
For almost all graphs $G$ in the Erd\"{o}s--R\'{e}nyi model $G_{n,p}$ with weight from the {\it Second multiple Zagreb} index $\ln \Pi_2(G)$,
$$\begin{array}{lll}
  \En(\ln \prod_2(G))&=&(1+o(1))\frac{16}{3\pi}\sqrt{p(1-p)}\cdot n^{3/2}\ln n.\\
  \end{array}$$
\end{cor}
It is easily seen that the maximum is achieved when $p=1/2$.\\

The Lanzhou index \cite{Lan} is defined as $Lz=(n-1)M_1(G) - F(G)$, where $M_1(G)$ is the first Zagreb index as defined
in (1), and $F(G)$ is called the forgotten index, which is defined as $F(G)=\sum_{ij\in E(G)} (d_i^2 +d_j^2)$.
Therefore,
$$Lz(G)=\sum_{ij\in E} [(n-1)(d_i+d_j)-(d_i^2 +d_j^2)]=\sum_{ij\in E}(d_i\bar{d_i} + d_j\bar{d_j}),$$
where $\bar{d_i}$ denotes the degree of vertex $i$ in the compliment $\bar{G}$ of $G$.
From the Lanzhou index, we can get the Lanzhou matrix as follows.\\

(12) For $LZ(G)=(LZ_{ij})$, $LZ_{ij}=(n-1)(d_i+d_j)-(d_i^2 +d_j^2)$ if $ij\in E$; 0 otherwise.
\begin{cor}\label{cor10}
For almost all graphs $G$ in the Erd\"{o}s--R\'{e}nyi model $G_{n,p}$ with weight from the Lanzhou index $Lz(G)$,
$$\begin{array}{lll}
  \En(LZ(G))=(1+o(1))\frac{16}{3\pi}\cdot[p(1-p)]^{3/2}\cdot n^{7/2}.\\
  \end{array}$$
\end{cor}
It is easily seen that the maximum is achieved when $p=1/2$.\\

\noindent {\bf Remark 2:} So far we can collect these matrices which are defined from indices of chemical use. For more
such indices we refer to \cite{FGD}.\\

\noindent {\bf Acknowledgement.} The authors are very grateful to the reviewers for their valuable comments and 
suggestions, which are very helpful to improving the paper.

\end{document}